# COMPONENT SELECTION AND SMOOTHING IN MULTIVARIATE NONPARAMETRIC REGRESSION

By Yi Lin[1] and Hao Helen Zhang[2]

*University of Wisconsin–Madison and North Carolina State University*

We propose a new method for model selection and model fitting in multivariate nonparametric regression models, in the framework of smoothing spline ANOVA. The "COSSO" is a method of regularization with the penalty functional being the sum of component norms, instead of the squared norm employed in the traditional smoothing spline method. The COSSO provides a unified framework for several recent proposals for model selection in linear models and smoothing spline ANOVA models. Theoretical properties, such as the existence and the rate of convergence of the COSSO estimator, are studied. In the special case of a tensor product design with periodic functions, a detailed analysis reveals that the COSSO does model selection by applying a novel soft thresholding type operation to the function components. We give an equivalent formulation of the COSSO estimator which leads naturally to an iterative algorithm. We compare the COSSO with MARS, a popular method that builds functional ANOVA models, in simulations and real examples. The COSSO method can be extended to classification problems and we compare its performance with those of a number of machine learning algorithms on real datasets. The COSSO gives very competitive performance in these studies.

**1. Introduction.** Consider the multivariate nonparametric regression problem $y_i = f(x_i) + \varepsilon_i, i = 1, \ldots, n$, where $f$ is the regression function to be estimated, $x_i = (x_i^{(1)}, \ldots, x_i^{(d)})$ are $d$-dimensional vectors of covariates and the $\varepsilon$ are independent noise variates with mean 0 and variance $\sigma^2$. The estimator is judged in terms of prediction accuracy and interpretability. A popular model for high-dimensional problems is the smoothing spline analysis of variance

Received April 2005; revised December 2005.
[1]Supported in part by NSF Grant DMS-01-34987.
[2]Supported in part by NSF Grant DMS-04-05913.

(SS-ANOVA) model [10, 19, 20]. In the SS-ANOVA we write

$$(1) \qquad f(x) = b + \sum_{j=1}^{d} f_j(x^{(j)}) + \sum_{j<k} f_{jk}(x^{(j)}, x^{(k)}) + \cdots,$$

where $b$ is a constant, the $f_j$'s are the main effects, the $f_{jk}$'s are the two-way interactions, and so on. The sequence is usually truncated somewhere to enhance interpretability. The identifiability of the terms in (1) is assured by side conditions through averaging operators. The SS-ANOVA generalizes the popular additive model and provides a general framework for nonparametric multivariate function estimation.

The common approach to estimation in SS-ANOVA is penalized least squares, with the penalty being a sum of squared norms of the terms in (1). One important question in the application of SS-ANOVA is to determine which variables or which ANOVA components should be included in the model. Gu [9] proposed using cosine diagnostics as model checking tools after model fitting in Gaussian regression. Chen [3] studied interaction spline models via SS-ANOVA and developed a nonstandard test procedure for model selection. Yau, Kohn and Wood [21] presented a Bayesian method for variable selection in a nonparametric manner. Gunn and Kandola [11] proposed a sparse kernel approach in a closely related framework. Zhang et al. [23] proposed a likelihood basis pursuit approach to model selection and estimation in the SS-ANOVA for exponential families. They expanded each nonparametric component function in (1) as a linear combination of a large number of basis functions and applied an $L_1$ penalty to the coefficients of all the basis functions. The $L_1$ penalty gives a solution that is sparse in the coefficients. However, a separate model selection procedure has to be applied after model fitting, since sparsity in coefficients helps but does not guarantee sparsity in SS-ANOVA components.

In this paper we introduce a new approach for model selection and estimation in SS-ANOVA. This is a penalized least squares method with the penalty functional being the sum of component norms, rather than the sum of squared component norms. Our method will be referred to as the COmponent Selection and Smoothing Operator (COSSO). The general methodology is introduced in Section 2, where we also prove the existence of the COSSO estimate and give some rate of convergence results. A connection between the COSSO and the popular LASSO in linear regression is shown. It turns out that when the COSSO formulation is used in linear models, it reduces to the LASSO. On the other hand, the COSSO gives an alternative interpretation of the penalty term in the LASSO to being the $L_1$ norm of the coefficients: it is the sum of component norms. Thus the COSSO can be seen as a nontrivial extension of the LASSO in linear models to multivariate nonparametric models. In Section 3 we obtain an alternative formulation of the COSSO that is more suitable for computation. In Section 4 we consider the special case of a tensor product design with periodic functions. A detailed analysis in this special case sheds light on the mechanism of the COSSO in terms of component selection in SS-ANOVA. In particular, we show in this case that the COSSO does model selection by applying a novel soft thresholding-type operation to the function components. In Section 5 we present a COSSO algorithm that is based on iterating between the smoothing spline method and the nonnegative garrote [1]. In Section 6 we consider the choice of the tuning parameter. Simulations are given in Section 7, where we compare the COSSO with the MARS procedure developed by Friedman [8], a popular algorithm that builds functional ANOVA models. The COSSO can be naturally extended to perform classification tasks, and we also compare the performance of the COSSO with that of many machine learning methods on some benchmark datasets. These real examples are given in Section 8, and Section 9 contains a discussion. The proofs are given in the Appendix.

## 2. The COSSO in smoothing spline ANOVA.

2.1. *The smoothing spline ANOVA.* In the commonly used smoothing spline ANOVA model over $\mathcal{X} = [0,1]^d$, it is assumed that $f \in \mathcal{F}$, where $\mathcal{F}$ is a reproducing kernel Hilbert space (RKHS) corresponding to the decomposition (1). Let $H^j$ be a function space of functions of $x^{(j)}$ over $[0,1]$ such that $H^j = \{1\} \oplus \bar{H}^j$. Then the tensor product space of the $H^j$'s is

$$(2) \qquad \bigotimes_{j=1}^{d} H^j = \{1\} \oplus \sum_{j=1}^{d} \bar{H}^j \oplus \sum_{j<k} [\bar{H}^j \otimes \bar{H}^k] \oplus \cdots.$$

Each functional component in the SS-ANOVA decomposition (1) lies in a subspace in the orthogonal decomposition (2) of $\bigotimes_{j=1}^{d} H^j$. Typically only low-order interactions are considered in the SS-ANOVA model for interpretability and visualization. The popular additive model is a special case in which $f(x^{(1)}, \ldots, x^{(d)}) = b + \sum_{j=1}^{d} f_j(x^{(j)})$, with $f_j \in \bar{H}^j$. In this case the selection of functional components is equivalent to variable selection. In more complex SS-ANOVA models model selection amounts to the selection of main effects and interaction terms in the SS-ANOVA decomposition. The interaction terms reside in the tensor product spaces of univariate function spaces. The reproducing kernel of a tensor product space is simply the product of the reproducing kernels of the individual spaces. This greatly facilitates the use of the smoothing spline-type method in such models.

A common example of the function space $H^j$ of univariate functions is the Sobolev Hilbert space. The $\ell$th-order Sobolev Hilbert space is: $S_\ell = \{g : g, g', \ldots, g^{(\ell-1)} \text{ are absolutely continuous}, g^{(\ell)} \in \mathcal{L}_2[0,1]\}$. Following [19], we define the norm in $S_\ell$ by

$$\|g\|^2 = \sum_{\nu=0}^{\ell-1} \left\{ \int_0^1 g^{(\nu)}(t)\,dt \right\}^2 + \int_0^1 \{g^{(\ell)}(t)\}^2\,dt. \tag{3}$$

With this norm $S_\ell$ can be decomposed as the direct sum of two orthogonal subspaces $S_\ell = \{1\} \oplus \bar{S}_\ell$. The spaces $S_\ell$ and $\bar{S}_\ell$ are RKHS's and their reproducing kernels are given in [19]. The second-order Sobolev Hilbert space $S_2$ is the most commonly used in practice and will be used in our implementation of the COSSO.

2.2. *The COSSO.* In general, the function space in the SS-ANOVA can be written as

$$\mathcal{F} = \{1\} \oplus \mathcal{F}_1 \qquad \text{with} \quad \mathcal{F}_1 = \bigoplus_{\alpha=1}^{p} \mathcal{F}^\alpha, \tag{4}$$

where $\mathcal{F}^1, \ldots, \mathcal{F}^p$ are $p$ orthogonal subspaces of $\mathcal{F}$. In the additive model $p = d$ and the $\mathcal{F}^\alpha$'s are the main effect spaces. In the two-way interaction model there are $d$ main effect spaces and $d(d-1)/2$ two-way interaction spaces, thus $p = d(d+1)/2$. We may further decompose the functional components into parametric and nonparametric parts, as is commonly done with the smoothing spline method. We do not pursue this in this paper as our emphasis is on the selection of functional components in SS-ANOVA. However, the general idea of our procedure can still be applied with this further decomposition, and it may be helpful to select parametric and nonparametric components of the variables.

Denote the norm in the RKHS $\mathcal{F}$ by $\|\cdot\|$. A traditional smoothing spline type method finds $f \in \mathcal{F}$ to minimize

$$\frac{1}{n}\sum_{i=1}^{n}\{y_i - f(x_i)\}^2 + \lambda \sum_{\alpha=1}^{p} \theta_\alpha^{-1} \|P^\alpha f\|^2, \tag{5}$$

where $P^\alpha f$ is the orthogonal projection of $f$ onto $\mathcal{F}^\alpha$ and $\theta_\alpha \geq 0$. If $\theta_\alpha = 0$, then the minimizer is taken to satisfy $\|P^\alpha f\|^2 = 0$. We use the convention $0/0 = 0$ throughout this paper. The smoothing parameter $\lambda$ is confounded with the $\theta$'s, but is usually included in the setup for computational purposes.

We propose the COSSO procedure that finds $f \in \mathcal{F}$ to minimize

$$\frac{1}{n}\sum_{i=1}^{n}\{y_i - f(x_i)\}^2 + \tau_n^2 J(f) \qquad \text{with } J(f) = \sum_{\alpha=1}^{p} \|P^\alpha f\|, \tag{6}$$

where $\tau_n$ is a smoothing parameter. We sometimes suppress the dependence of $\tau$ on $n$ in our notation. The penalty term $J(f)$ in the COSSO is a sum of RKHS norms, instead of the squared RKHS norm penalty employed in the smoothing spline. The penalty $J(f)$ is not a norm in $\mathcal{F}$. However, it is a convex functional and is a pseudonorm in the sense: for any $f, g$ in $\mathcal{F}, J(f) \geq 0, J(cf) = |c|J(f), J(f+g) \leq J(f) + J(g)$; for any nonconstant $f$ in $\mathcal{F}, J(f) > 0$. And we have that

$$\text{(7)} \qquad \sum_{\alpha=1}^{p} \|P^\alpha f\|^2 \leq J^2(f) \leq p \sum_{\alpha=1}^{p} \|P^\alpha f\|^2.$$

The existence of the COSSO estimate is guaranteed due to the convexity of (6), as stated in the following.

THEOREM 1. *Let $\mathcal{F}$ be an RKHS of functions over an input space $\mathcal{X}$. Assume that $\mathcal{F}$ can be decomposed as in (4). Then there exists a minimizer of (6) in $\mathcal{F}$.*

2.3. *Connection to the LASSO in linear models.* In linear models the regression function is assumed to be $f(x) = \beta_0 + \sum_{j=1}^{d} \beta_j x^{(j)}$. Traditional approaches to variable selection include best subset selection and forward/backward stepwise selection. As pointed out by [1], these methods suffer from instability and relative lack of accuracy. Several new and effective methods for variable selection in linear models have been proposed in recent years [1, 5, 6, 7, 15]. Two methods, the nonnegative garrote [1] and the LASSO [15], are related to the method in our paper, and are reviewed in the following.

Let $\hat{\beta}^o = (\hat{\beta}_0^o, \ldots, \hat{\beta}_d^o)$ be the ordinary least squares estimates. The nonnegative garrote solution is $(\hat{\beta}_0^o, r_1 \hat{\beta}_1^o, \ldots, r_d \hat{\beta}_d^o)$, where $(r_1, \ldots, r_d)$ is the solution to

$$\min_{r_1, \ldots, r_d} \sum_{i=1}^{n} \left\{ y_i - \hat{\beta}_0^o - \sum_{j=1}^{d} r_j \hat{\beta}_j^o x_i^{(j)} \right\}^2$$

$$\text{subject to} \quad r_j \geq 0, j = 1, \ldots, d \text{ and } \sum_{j=1}^{d} r_j \leq t.$$

Here $t \geq 0$ is a tuning parameter. The nonnegative garrote selects a subset and shrinks the estimate at the same time. Breiman [1] showed that the nonnegative garrote has consistently lower prediction error than subset selection with extensive simulation studies.

The Least Absolute Shrinkage and Selection Operator (LASSO) estimate $\hat{\beta} = (\hat{\beta}_0, \ldots, \hat{\beta}_d)$ is the minimizer of

$$\frac{1}{n} \sum_{i=1}^{n} \left\{ y_i - \beta_0 - \sum_{j=1}^{d} \beta_j x_i^{(j)} \right\}^2 \quad \text{subject to} \quad \sum_{j=1}^{d} |\beta_j| \leq t,$$

or equivalently, the minimizer of

$$\frac{1}{n}\sum_{i=1}^{n}\left\{y_i - \beta_0 - \sum_{j=1}^{d}\beta_j x_i^{(j)}\right\}^2 + \lambda\sum_{j=1}^{d}|\beta_j|,$$

where $t$ and $\lambda$ are tuning parameters. The LASSO is a penalized least squares method with the $L_1$ penalty on the coefficients.

The LASSO can be seen as a special case of the COSSO. For the input space $\mathcal{X} = [0,1]^d$, consider the linear function space $\mathcal{F} = \{1\} \oplus \{x^{(1)} - 1/2\} \oplus \cdots \oplus \{x^{(d)} - 1/2\}$, with the usual $L_2$ inner product on $\mathcal{F}$: $(f,g) = \int_{\mathcal{X}} fg$. The penalty term in the COSSO (6) becomes $J(f) = (12)^{-1/2}\sum_{j=1}^{d}|\beta_j|$ for $f(x) = \beta_0 + \sum_{j=1}^{d}\beta_j x^{(j)}$. This is equivalent to the $L_1$ norm on the linear coefficients, leading to the LASSO estimator. Notice, however, in the COSSO interpretation the penalty is the sum of the norms of the function components, rather than the $L_1$ norm of the coefficients.

2.4. *Asymptotic property of the COSSO.* In this section we assume a fixed design. Define $y = (y_1, \ldots, y_n)^{\mathrm{T}}$. With a little abuse of notation, let $f$ stand for both the regression function and its functional values at data points, that is, $f = (f(x_1), \ldots, f(x_n))^{\mathrm{T}}$. Define the norm $\|\cdot\|_n$ and inner product $\langle \cdot, \cdot \rangle_n$ in $R^n$ as

$$\|f\|_n^2 = \frac{1}{n}\sum_{i=1}^{n} f^2(x_i), \qquad \langle f, g \rangle_n = \frac{1}{n}\sum_{i=1}^{n} f(x_i)g(x_i);$$

then $\|y - f\|_n^2 = 1/n \sum_{i=1}^{n} \{y_i - f(x_i)\}^2$. The following theorem shows that the COSSO estimator in the additive model has a rate of convergence $n^{-\ell/(2\ell+1)}$, where $\ell$ is the order of smoothness of the components.

THEOREM 2. *Consider the regression model $y_i = f_0(x_i) + \varepsilon_i$, $i = 1, \ldots, n$, where the $x_i$'s are given covariates in $[0,1]^d$, and the $\varepsilon_i$'s are independent $N(0, \sigma^2)$ noise. Assume $f_0$ lies in $\mathcal{F} = \{1\} \oplus \mathcal{F}_1$, $\mathcal{F}_1 = \bigoplus_{j=1}^{d} \bar{S}^j$, with $S^j = \{1\} \oplus \bar{S}^j$ being the $\ell$th-order Sobolev space with norm (3). Consider the COSSO estimate $\hat{f}$ as defined in (6). Then (i) if $f_0$ is not a constant, and $\tau_n^{-1} = O_p(n^{\ell/(2\ell+1)})J^{(2\ell-1)/(4\ell+2)}(f_0)$, we have $\|\hat{f} - f_0\|_n = O_p(\tau_n)J^{1/2}(f_0)$; (ii) if $f_0$ is a constant, we have $\|\hat{f} - f_0\|_n = O_p(\max\{n^{-\ell/(2\ell-1)}\tau_n^{-2/(2\ell-1)}, n^{-1/2}\})$.*

**3. An equivalent formulation.** It can be shown that the solution to (6) is in a finite-dimensional space, and therefore the COSSO estimate can be computed directly from (6).

LEMMA 1. *Let $\hat{f} = \hat{b} + \sum_{\alpha=1}^{p} \hat{f}_\alpha$ be a minimizer of (6) in (4), with $\hat{f}_\alpha \in \mathcal{F}^\alpha$. Then $\hat{f}_\alpha \in \mathrm{span}\{R_\alpha(x_i, \cdot), i = 1, \ldots, n\}$, where $R_\alpha(\cdot, \cdot)$ is the reproducing kernel of $\mathcal{F}^\alpha$.*

However, it is possible to give an equivalent form of (6) that is easier to compute. Consider the problem of finding $\theta = (\theta_1, \ldots, \theta_p)^{\mathrm{T}}$ and $f \in \mathcal{F}$ to minimize

$$\frac{1}{n}\sum_{i=1}^{n}\{y_i - f(x_i)\}^2 + \lambda_0 \sum_{\alpha=1}^{p} \theta_\alpha^{-1}\|P^\alpha f\|^2 + \lambda \sum_{\alpha=1}^{p} \theta_\alpha \quad \text{subject to} \quad \theta_\alpha \geq 0,$$
(8)
$$\alpha = 1, \ldots, p,$$

where $\lambda_0$ is a constant and $\lambda$ is a smoothing parameter. The constant $\lambda_0$ can be fixed at any positive value and is included here only for computational considerations.

LEMMA 2. *Set $\lambda = \tau^4/(4\lambda_0)$. (i) If $\hat{f}$ minimizes (6), set $\hat{\theta}_\alpha = \lambda_0^{1/2}\lambda^{-1/2} \times \|P^\alpha \hat{f}\|$; then the pair $(\hat{\theta}, \hat{f})$ minimizes (8). (ii) On the other hand, if a pair $(\hat{\theta}, \hat{f})$ minimizes (8), then $\hat{f}$ minimizes (6).*

The form of (8) is very similar to the common smoothing spline (5) with multiple smoothing parameters, except that there is an additional penalty on the $\theta$'s. Differently from [23] where the sparsity in coefficients does not assure the sparsity in the functional components, for the COSSO procedure the sparsity of each component $f_j$ is controlled by a single parameter $\theta_j$. The additional penalty term on the $\theta$'s, $\sum_{j=1}^{p} \theta_j$, shrinks them toward zero and hence makes the $\theta$'s sparse, giving rise to zero function components in the COSSO estimate.

We only need to tune $\lambda$ in the implementation of the COSSO, and $\lambda_0$ can be fixed at any positive value. In Section 6 we will see that an appropriate choice of $\lambda_0$ helps to put the tuning parameter on a natural scale and facilitates the tuning of the COSSO. In contrast, the common smoothing spline has two sets of smoothing parameters $\lambda$ and $\theta$'s that are confounded. The common way to search for smoothing parameters iterates between $\lambda$ and the $\log\theta$'s, making it difficult to have zero components in the solution.

**4. A special case with the tensor product design.** To illustrate the mechanism of the COSSO for model selection, in this section we give an instructive analysis in the special case of a tensor product design with an SS-ANOVA model built from the second-order Sobolev spaces of periodic functions. We assume the $\varepsilon$'s in the regression model are independent with distribution $N(0, \sigma^2)$. In a tensor product design the design points are

$$\{(x_{i_1,1}, x_{i_2,2}, \ldots, x_{i_d,d}) : i_k = 1, \ldots, n_k, k = 1, \ldots, d\},$$

where $x_{j,k} = j/n_k$, $j = 1, \ldots, n_k$, $k = 1, \ldots, d$. Without loss of generality, we fix $\lambda_0 = 1$ in the COSSO (8) and focus on the case $d = 2$ with the SS-ANOVA model being $f(s,t) = b + f_1(s) + f_2(t) + f_{12}(s,t)$. We assume $n_1 = n_2 = m$ is an even integer. The sample size is then $n = m^2$.

The second-order Sobolev space of periodic functions can be written as $T = \{1\} \oplus \bar{T}$, where

$$\bar{T} = \left\{ f : f(t) = \sum_{\nu=1}^{\infty} a_\nu \sqrt{2} \cos 2\pi\nu t + \sum_{\nu=1}^{\infty} b_\nu \sqrt{2} \sin 2\pi\nu t, \right.$$
$$\left. \text{with } \sum_{\nu=1}^{\infty} (a_\nu^2 + b_\nu^2)(2\pi\nu)^4 < \infty \right\}.$$

The norm in $\bar{T}$ is $\|g\|^2 = \int_0^1 \{g''(t)\}^2 \, dt$. When $m$ is large, a good approximate subspace of $T$ is $T_m = \{1\} \oplus \bar{T}_m$ with

$$\bar{T}_m = \left\{ f : f(t) = \sum_{\nu=1}^{m/2-1} a_\nu \sqrt{2} \cos 2\pi\nu t + \sum_{\nu=1}^{m/2-1} b_\nu \sqrt{2} \sin 2\pi\nu t + a_{m/2} \cos \pi m t \right\}.$$

Wahba [19] used this subspace approximation to give a very instructive investigation of the filtering properties of the smoothing spline. Here we consider minimizing (8) in $\mathcal{F}_m = T_m^1 \otimes T_m^2 = \{1\} \oplus \bar{T}_m^1 \oplus \bar{T}_m^2 \oplus (\bar{T}_m^1 \otimes \bar{T}_m^2)$, as the argument is instructive. The argument for the more general function space $\mathcal{F} = T^1 \otimes T^2$ is similar but involves more technicality, and is deferred to the Appendix.

In this case (8) becomes

$$\frac{1}{n} \sum_{k=1}^{m} \sum_{\ell=1}^{m} \{y_{kl} - f(x_{k,1}, x_{\ell,2})\}^2$$
$$+ \theta_1^{-1} \int_0^1 \left\{ \frac{\partial^2 f_1(s)}{\partial s^2} \right\}^2 ds + \theta_2^{-1} \int_0^1 \left\{ \frac{\partial^2 f_2(t)}{\partial t^2} \right\}^2 dt$$
$$+ \theta_{12}^{-1} \int_0^1 \int_0^1 \left\{ \frac{\partial^4 f_{12}(s,t)}{\partial s^2 \, \partial t^2} \right\}^2 ds \, dt + \lambda(\theta_1 + \theta_2 + \theta_{12}),$$
$$\theta_1 \geq 0, \theta_2 \geq 0, \theta_{12} \geq 0.$$

Write $\gamma_1(t) = 1$, $\gamma_{2\nu}(t) = \sqrt{2} \cos(2\pi\nu t)$, $\gamma_{2\nu+1}(t) = \sqrt{2} \sin(2\pi\nu t)$ for $\nu = 1, \ldots, m/2 - 1$, and $\gamma_m(t) = \cos(\pi m t)$. Then any function in $T_m$ can be written as $g(t) = \sum_{\nu=1}^{m} a_\nu \gamma_\nu(t)$, and any function in $\mathcal{F}_m$ can be written as

(9)  $$f(s,t) = \sum_{\mu=1}^{m} \sum_{\nu=1}^{m} a_{\mu\nu} \gamma_\mu(s) \gamma_\nu(t).$$

It is known that (see [19], page 23)

$$m^{-1} \sum_{k=1}^{m} \gamma_\mu(k/m) \gamma_\nu(k/m) = \begin{cases} 1, & \text{if } \mu = \nu = 1, \ldots, m, \\ 0, & \text{if } \mu \neq \nu, \mu, \nu = 1, \ldots, m. \end{cases}$$

Recall the inner product $\langle \cdot, \cdot \rangle_n$ of $R^n$ defined in Section 2.4. Write $\gamma_{\mu\nu}(s,t) = \gamma_\mu(s)\gamma_\nu(t)$, and $\gamma_{\mu\nu}$ as the data vector corresponding to the function $\gamma_{\mu\nu}(s,t)$. From the above orthogonality relations and the tensor product design, we get

$$\langle \gamma_{\mu_1\nu_1}, \gamma_{\mu_2\nu_2} \rangle_n = \begin{cases} 1, & \text{if } \mu_1 = \mu_2 = 1, \ldots, m; \nu_1 = \nu_2 = 1, \ldots, m, \\ 0, & \text{if } \mu_1 \neq \mu_2 \text{ or } \nu_1 \neq \nu_2, \mu_1, \nu_1, \mu_2, \nu_2 = 1, \ldots, m. \end{cases}$$

Therefore $\{\gamma_{\mu\nu}, \mu = 1, \ldots, m; \nu = 1, \ldots, m\}$ form an orthonormal basis in $R^n$ with respect to the norm $\|\cdot\|_n$. We then get from (9) that $a_{\mu\nu} = \langle f, \gamma_{\mu\nu} \rangle_n$. Write $z_{\mu\nu} = \langle y, \gamma_{\mu\nu} \rangle_n$. Then $z_{\mu\nu} = a_{\mu\nu} + \delta_{\mu\nu}$, where $\delta_{\mu\nu} \sim N(0, \sigma^2/n)$ are independent. The COSSO problem can be written as

$$(10) \quad \sum_{\mu=1}^{m}\sum_{\nu=1}^{m}(z_{\mu\nu} - a_{\mu\nu})^2 + \theta_1^{-1}\sum_{\mu=2}^{m}q_{\mu 1}a_{\mu 1}^2 + \theta_2^{-1}\sum_{\nu=2}^{m}q_{1\nu}a_{1\nu}^2$$
$$+ \theta_{12}^{-1}\sum_{\mu=2}^{m}\sum_{\nu=2}^{m}q_{\mu\nu}a_{\mu\nu}^2 + \lambda(\theta_1 + \theta_2 + \theta_{12}),$$

with $q_{\mu\nu} \sim \mu^4\nu^4$ uniformly for $\mu \neq 1$ or $\nu \neq 1$, $\mu, \nu = 1, \ldots, m$. Here $\sim$ is read as "has the same order as." Therefore the minimizing $a_{\mu\nu}$ satisfies $\hat{a}_{11} = z_{11}$; $\hat{a}_{\mu 1} = z_{\mu 1}\theta_1(\theta_1 + q_{\mu 1})^{-1}$, for $\mu \geq 2$; $\hat{a}_{1\nu} = z_{1\nu}\theta_2(\theta_2 + q_{1\nu})^{-1}$, for $\nu \geq 2$; $\hat{a}_{\mu\nu} = z_{\mu\nu}\theta_{12}(\theta_{12} + q_{\mu\nu})^{-1}$, for $\mu \geq 2, \nu \geq 2$; and (10) becomes

$$\left\{\sum_{\mu=2}^{m} q_{\mu 1}z_{\mu 1}^2(q_{\mu 1} + \theta_1)^{-1} + \lambda\theta_1\right\} + \left\{\sum_{\nu=2}^{m} q_{1\nu}z_{1\nu}^2(q_{1\nu} + \theta_2)^{-1} + \lambda\theta_2\right\}$$
$$+ \left\{\sum_{\mu=2}^{m}\sum_{\nu=2}^{m} q_{\mu\nu}z_{\mu\nu}^2(q_{\mu\nu} + \theta_{12})^{-1} + \lambda\theta_{12}\right\}.$$

We see that the three components can be minimized separately. Let us concentrate on $\theta_{12}$, as $\theta_1$ and $\theta_2$ can be dealt with similarly. Let

$$A(\theta_{12}) = \sum_{\mu=2}^{m}\sum_{\nu=2}^{m} q_{\mu\nu}z_{\mu\nu}^2(q_{\mu\nu} + \theta_{12})^{-1} + \lambda\theta_{12}.$$

Then $A'(\theta_{12}) = \lambda - \sum_{\mu=2}^{m}\sum_{\nu=2}^{m} q_{\mu\nu}z_{\mu\nu}^2(q_{\mu\nu} + \theta_{12})^{-2}$, which increases as $\theta_{12} \geq 0$ increases. Define $U = \sum_{\mu=2}^{m}\sum_{\nu=2}^{m} q_{\mu\nu}^{-1}z_{\mu\nu}^2$. If $U \leq \lambda$, then $A'(0) \geq 0$, $A'(\theta_{12}) > 0$ for all $\theta_{12} > 0$, and the minimizing $\hat{\theta}_{12}$ of $A$ is 0; otherwise the minimizing $\hat{\theta}_{12}$ is larger than 0. Therefore we can see that the COSSO estimator selects components through a soft thresholding-type operation according to the magnitude of $U$. Notice $\hat{\theta}_{12} = 0$ implies $\hat{f}_{12} = 0$.

With the analysis above, we can now show that, when $\lambda \to 0$ and $n\lambda \to \infty$, the COSSO selects the correct model with probability tending to 1. Without loss of generality, let us concentrate on $f_{12}$.

If $f_{12} = 0$, then $a_{\mu\nu} = 0$ for any pair $(\mu,\nu)$ such that $\mu \geq 2$ and $\nu \geq 2$. So $E(U) \sim \sigma^2/n \sum_{\mu=2}^{m} \sum_{\nu=2}^{m} \mu^{-4}\nu^{-4} \sim n^{-1}\sigma^2$, $\mathrm{var}(U) = \sum_{\mu=2}^{m}\sum_{\nu=2}^{m} 2n^{-2} \times \sigma^4 q_{\mu\nu}^{-2} \sim n^{-2}\sigma^4$. Therefore when $n\lambda \to \infty$, by Chebyshev's inequality,

$$\mathrm{pr}(U > \lambda) \leq \mathrm{pr}(|U - E(U)| > \lambda - E(U)) \leq \mathrm{var}(U)/\{\lambda - E(U)\}^2 \to 0.$$

Therefore with probability tending to unity, $U \leq \lambda$, and thus $\hat{f}_{12} = 0$.

On the other hand, if $f_{12} \neq 0$, then $a_{\mu_0,\nu_0} \neq 0$ for some $\mu_0 \geq 2$ and $\nu_0 \geq 2$, and

$$E(U) \geq E(q_{\mu_0,\nu_0}^{-1} z_{\mu_0,\nu_0}^2) \geq q_{\mu_0,\nu_0}^{-1} a_{\mu_0,\nu_0}^2,$$

$$\mathrm{var}(U) = \sum_{\mu=2}^{m}\sum_{\nu=2}^{m} q_{\mu\nu}^{-2} \mathrm{var}(z_{\mu\nu}^2) = \sum_{\mu=2}^{m}\sum_{\nu=2}^{m} q_{\mu\nu}^{-2} (4n^{-1} a_{\mu\nu}^2 \sigma^2 + 2n^{-2}\sigma^4)$$

$$\leq \left\{ 4n^{-1}\sigma^2 \sum_{\mu=2}^{m}\sum_{\nu=2}^{m} a_{\mu\nu}^2 \right\} + 2n^{-2}\sigma^4$$

$$= 4n^{-1}\sigma^2 \|f_{12}\|_{L_2}^2 + 2n^{-2}\sigma^4 = O(n^{-1}).$$

Therefore when $\lambda \to 0$, by Chebyshev's inequality, we get

$$\mathrm{pr}(U < \lambda) \leq \mathrm{pr}(|U - E(U)| > E(U) - \lambda) \leq \mathrm{var}(U)/\{E(U) - \lambda\}^2 \to 0.$$

Therefore with probability tending to unity, $U > \lambda$, and thus $\hat{f}_{12} \neq 0$.

**5. Algorithm.** For any fixed $\theta$, the COSSO (8) is equivalent to the smoothing spline (5). Therefore from the smoothing spline literature (e.g., [19]) it is well known the solution $f$ has the form $f(x) = \sum_{i=1}^{n} c_i R_\theta(x_i, x) + b$, where $c = (c_1, \ldots, c_n)^{\mathrm{T}} \in R^n$, $b \in R$, and $R_\theta = \sum_{\alpha=1}^{p} \theta_\alpha R_\alpha$, with $R_\alpha$ being the reproducing kernel of $\mathcal{F}^\alpha$. With some abuse of notation, let $R_\alpha$ also stand for the $n \times n$ matrix $\{R_\alpha(x_i, x_j)\}$, $i = 1, \ldots, n$, $j = 1, \ldots, n$, let $R_\theta$ also stand for the matrix $\sum_{\alpha=1}^{p} \theta_\alpha R_\alpha$, and let $\mathbf{1}_r$ be the column vector consisting of $r$ 1's. Then we can write $f = R_\theta c + b\mathbf{1}_n$, and (8) can be expressed as

$$(11) \quad \frac{1}{n}\left(y - \sum_{\alpha=1}^{p} \theta_\alpha R_\alpha c - b\mathbf{1}_n\right)^{\mathrm{T}} \left(y - \sum_{\alpha=1}^{p} \theta_\alpha R_\alpha c - b\mathbf{1}_n\right) + \lambda_0 \sum_{\alpha=1}^{p} \theta_\alpha c^{\mathrm{T}} R_\alpha c + \lambda \sum_{\alpha=1}^{p} \theta_\alpha,$$

where $\theta_\alpha \geq 0$, $\alpha = 1, \ldots, p$.

If the $\theta$'s are fixed, then (11) can be written as

$$(12) \quad \min_{c,b} (y - R_\theta c - b\mathbf{1}_n)^{\mathrm{T}} (y - R_\theta c - b\mathbf{1}_n) + n\lambda_0 c^{\mathrm{T}} R_\theta c.$$

The solution to this smoothing spline problem is given in [19].

On the other hand, if $c$ and $b$ are fixed, denote $g_\alpha = R_\alpha c$ and let $G$ be the $n \times p$ matrix with the $\alpha$th column being $g_\alpha$. Simple calculation shows that the $\theta = (\theta_1, \ldots, \theta_p)^\mathrm{T}$ that minimizes (11) is the solution to

$$(13) \quad \min_\theta (z - G\theta)^\mathrm{T}(z - G\theta) + n\lambda \sum_{\alpha=1}^p \theta_\alpha \quad \text{subject to} \quad \theta_\alpha \geq 0, \alpha = 1, \ldots, p,$$

where $z = y - (1/2)n\lambda_0 c - b\mathbf{1}_n$.

Therefore a reasonable scheme would be to iterate between (12) and (13). In each iteration (11) is decreased. Notice that (13) is equivalent to

$$(14) \quad \min_\theta (z - G\theta)^\mathrm{T}(z - G\theta) \quad \text{subject to} \quad \theta_\alpha \geq 0, \alpha = 1, \ldots, p; \sum_{\alpha=1}^p \theta_\alpha \leq M,$$

for some $M \geq 0$. We prefer to iterate between (12) and (14) for computational considerations.

Notice that the formulation (14) is exactly the problem in calculating the nonnegative garrote estimate. Therefore our algorithm iterates between the smoothing spline and the nonnegative garrote. The algorithm starts with a natural initial solution given by the smoothing spline, which is already a good estimate. By applying later iterations of our algorithm, we get what we view as an iterative improvement on the smoothing spline. A limited number of iterations is usually sufficient to achieve good performance in practical applications. This is in spirit similar to the basis pursuit algorithm in [2]. We observe empirically that the COSSO objective function decreases quickly in the first iteration, and the objective function after the first iteration is already very close to the objective function at convergence, as the magnitude of the decrease in the first iteration dominates the decreases in subsequent iterations. This motivates us to consider the following one-step update procedure:

1. Initialization: Fix $\theta_\alpha = 1$, $\alpha = 1, \ldots, p$.
2. Solve for $c$ and $b$ with (12).
3. For $c$ and $b$ obtained in step 2, solve for $\theta$ with the nonnegative garrote (14).
4. With the new $\theta$, solve for $c$ and $b$ with the smoothing spline (12).

This one-step update procedure has the flavor of the one-step maximum likelihood procedure in which a one-step Newton–Raphson algorithm is applied to a good initial estimator and which is as efficient as fully iterated maximum likelihood. A discussion of the one-step procedure and the fully iterated procedure (in a different algorithm) can be found in [6]. In our experience, the one-step update procedure and the fully iterated procedure have comparable estimation accuracy.

**6. Choosing the tuning parameter.** The generalized cross-validation proposed by Craven and Wahba [4] is one of the most popular methods for choosing smoothing parameters in the smoothing spline method. Let $A$ be the smoothing matrix of the smoothing spline. That is, $\hat{y} = Ay$. The generalized cross-validation estimate of the risk is

$$GCV = \frac{\|\hat{y} - y\|_n^2}{\{n^{-1}\operatorname{tr}(I-A)\}^2}.$$

Tibshirani [15] proposed a GCV-type criterion for choosing the tuning parameter for the LASSO through a ridge estimate approximation. This approximation is particularly easy to understand in light of the form (8) for the linear model $f(x) = \beta_0 + \sum_{j=1}^d \beta_j x^{(j)}$: fix the $\theta_j$ at their estimated values $\hat{\theta}_j$, and calculate GCV for the corresponding ridge regression. This approximation ignores some variability in the estimation process. However, the simulation study in [15] suggests that it is a useful approximation. This motivates our GCV-type criterion: We use the GCV score for the smoothing spline in (8) when the $\theta$'s are fixed at the solution.

Another popular technique for choosing tuning parameters is fivefold or tenfold cross-validation. The computation load of GCV is smaller. We compare the performance of these two criteria in the COSSO with simulations. It is also possible to use the $C_p$ criterion based on the concept of generalized degrees of freedom [13, 22]. We do not consider this possibility since in our problem there is no explicit formula for the degrees of freedom and numerical evaluations tend to be computationally intensive.

The following is the complete algorithm for the COSSO with adaptive tuning:

1. Fix $\theta_\alpha = 1$, $\alpha = 1, \ldots, p$. Solve the smoothing spline problem, and tune $\lambda_0$ according to CV or GCV. Fix $\lambda_0$ at the chosen value in all later steps.
2. For each fixed $M$ in a reasonable range, apply the one-step COSSO algorithm with $M$. Choose the best $M$ according to CV or GCV. The solution corresponding to this chosen $M$ is the final solution.

In our simulations it is noticed that once $\lambda_0$ is fixed according to step 1, the optimal $M$ seems to be close to the number of important components. This helps to determine the range of tuning for $M$.

**7. Simulations.** In this section we study the empirical performance of the COSSO estimate in terms of estimation accuracy and model selection. We compare the COSSO with GCV, with fivefold cross-validation, and with MARS, which is a popular stepwise forward–backward procedure for building functional ANOVA models. The measure of accuracy is the integrated squared error ISE $= E_X\{\hat{f}(X) - f(X)\}^2$, which is estimated by Monte Carlo integration using 10,000 test points from the same distribution as the training points. We run each simulation example 100 times and average. The matlab code for the COSSO is available from webpages of the authors (www.stat.wisc.edu/~yilin or www4.stat.ncsu.edu/~hzhang). The MARS simulations are done in R, with the function "mars" in the "mda" library contributed by Trevor Hastie and Robert Tibshirani.

The following four functions on $[0, 1]$ are used as building blocks of regression functions in some of the simulations: $g_1(t) = t$; $g_2(t) = (2t - 1)^2$; $g_3(t) = \frac{\sin(2\pi t)}{2-\sin(2\pi t)}$; and $g_4(t) = 0.1\sin(2\pi t) + 0.2\cos(2\pi t) + 0.3\sin^2(2\pi t) + 0.4\cos^3(2\pi t) + 0.5\sin^3(2\pi t)$. Consider two covariance structures of the input vector $X$, with varying correlation:

*Compound symmetry.* Let $X^{(j)} = (W_j + tU)/(1 + t)$, $j = 1, \ldots, d$, where $W_1, \ldots, W_d$ and $U$ are i.i.d. from Uniform$(0, 1)$. Therefore corr$(X^{(j)}, X^{(k)}) = t^2/(1 + t^2)$ for $j \neq k$. The uniform design corresponds to the case $t = 0$.

(*Trimmed*) AR(1). Let $W_1, \ldots, W_d$ be i.i.d. $N(0, 1)$, and $X^{(1)} = W_1$, $X^{(j)} = \rho X^{(j-1)} + (1 - \rho^2)^{1/2} W_j$, $j = 2, \ldots, d$. Trim $X^{(j)}$ in $[-2.5, 2.5]$ and scale to $[0, 1]$.

EXAMPLE 1. Consider a simple additive model in $R^{10}$, with the underlying regression function $f(x) = 5g_1(x^{(1)}) + 3g_2(x^{(2)}) + 4g_3(x^{(3)}) + 6g_4(x^{(4)})$. Therefore $X^{(5)}, \ldots, X^{(10)}$ are uninformative. We consider the sample size $n = 100$. Generate $y = f(x) + \varepsilon$, where $\varepsilon$ is distributed as $N(0, 1.74)$. The standard deviation of the noise was chosen to give a signal-to-noise ratio 3:1 in the uniform case. For comparison, the variances of the component functions are var$\{5g_1(X^{(1)})\} = 2.08$, var$\{3g_2(X^{(2)})\} = 0.80$, var$\{4g_3(X^{(3)})\} = 3.30$ and var$\{6g_4(X^{(4)})\} = 9.45$.

We apply the COSSO with additive models (the additive COSSO) to the simulated data. Therefore there are ten functional components in the model. Figure 1 shows how the magnitudes of the estimated components change with the tuning parameter $M$ in one run. The magnitudes of the functional components are measured by their empirical $L_1$ norms, defined as $1/n \sum_{i=1}^{n} |\hat{f}_j(x_i^{(j)})|$ for $j = 1, \ldots, d$. The $\lambda_0$ in this run is fixed at $9.7656 \times 10^{-6}$. Both GCV and fivefold cross-validation choose $M = 3.5$, giving a model of five terms in this run. The estimated function components are plotted along with the true function components in Figure 2. Notice the components are centered according to the ANOVA decomposition.

For each setting of covariance structure, we run the simulation 100 times and average. The resulting average integrated squared error and its associated standard error (in parentheses) are given in Table 1. Also included in the table is the average integrated squared error of MARS for additive models. We can see the two COSSO procedures perform better than MARS in all the settings studied. To study the performance of the COSSO in terms of model selection, we determine in the uniform case the number of times each variable appears in the 100 chosen models (Table 2), and the number of terms in the 100 chosen models (Table 3). In our calculation we take $\theta$ to be zero if it is smaller than $10^{-6}$. The COSSO with fivefold cross-validation misses the second variable six times, but chooses the correct four-variable model 84 times. The COSSO with GCV and MARS do not miss any important variable, but tend to include uninformative variables in the chosen models. The COSSO with GCV chooses the correct four-variable model 57 times, while MARS does so only four times.

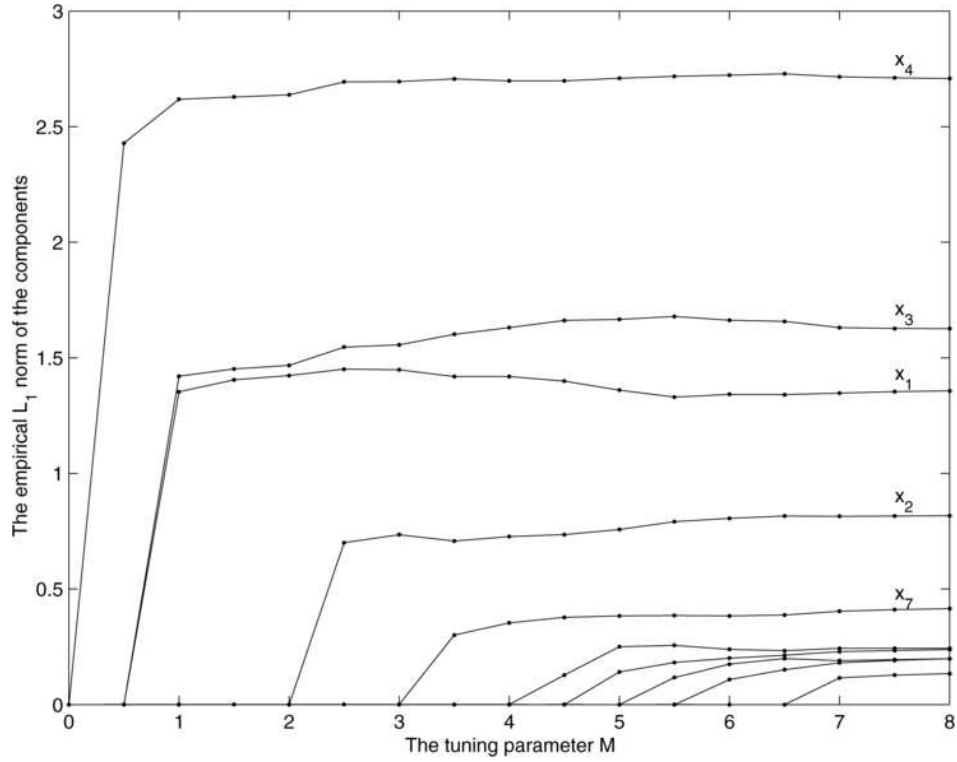

FIG. 1. *The empirical $L_1$ norm of the estimated components as plotted against the tuning parameter $M$ in one run of Example* 1.

TABLE 1
*Comparison of average integrated squared errors for Example* 1

|  | Comp. symm. | | | AR(1) | | |
| --- | --- | --- | --- | --- | --- | --- |
|  | $t=0$ | $t=1$ | $t=3$ | $\rho=-0.5$ | $\rho=0$ | $\rho=0.5$ |
| COSSO(GCV) | 0.93 (0.05) | 0.92 (0.04) | 0.97 (0.07) | 0.94 (0.05) | 1.04 (0.07) | 0.98 (0.07) |
| COSSO(5CV) | 0.80 (0.03) | 0.97 (0.05) | 1.07 (0.06) | 1.03 (0.06) | 1.03 (0.06) | 0.98 (0.05) |
| MARS | 1.57 (0.07) | 1.24 (0.06) | 1.30 (0.06) | 1.32 (0.07) | 1.34 (0.07) | 1.36 (0.08) |

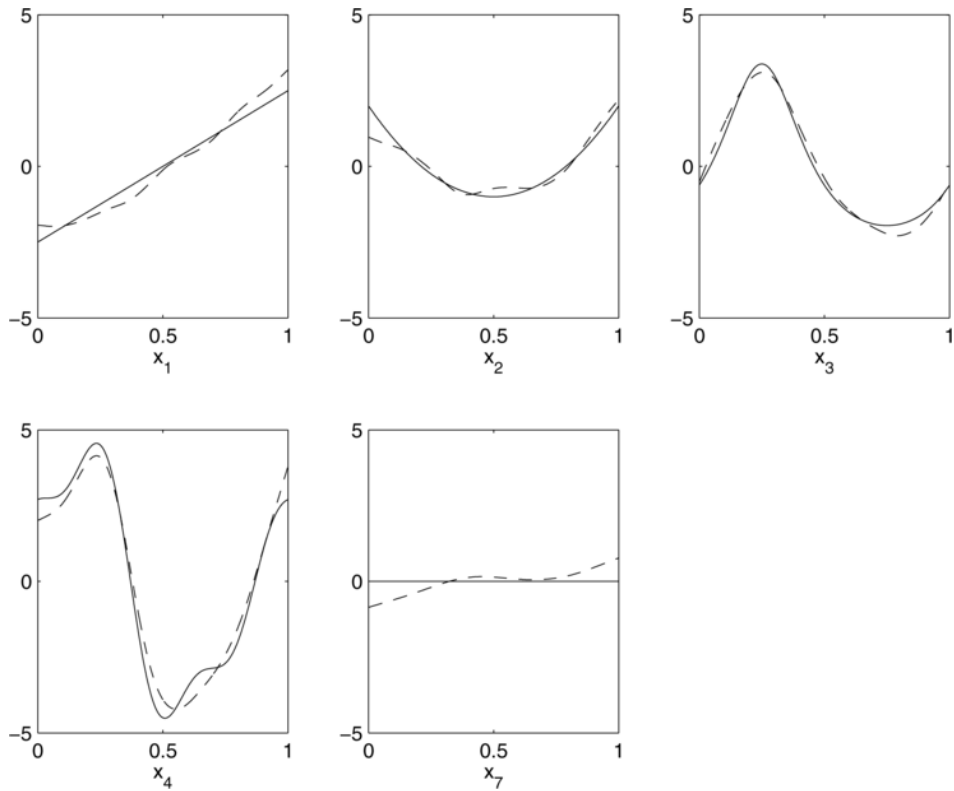

Fig. 2. *The estimated component (dashed line) and true component (solid line) functions in one run of Example* 1. *Shown are the components for variables* 1, 2, 3, 4 *and* 7. *For the other variables, both true and estimated components are zero.*

Table 4 gives the mean and standard deviation of the model sizes chosen by the methods in various settings. The settings considered are compound symmetry with $t=1$ and 3 and trimmed AR(1) with $\rho = -0.5$, 0 and 0.5. The average model size chosen by the COSSO with fivefold cross-validation is close to 4, the size of the true model. The COSSO with GCV selects slightly larger models. The models chosen by MARS are even larger.

EXAMPLE 2. Consider a larger model with $d = 60$ and the regression function
$$f(x) = g_1(x^{(1)}) + g_2(x^{(2)}) + g_3(x^{(3)}) + g_4(x^{(4)})$$

TABLE 2
*Appearance frequency of the variables in the models in the uniform setting*

|  | Variable | | | | | | | | | |
| --- | --- | --- | --- | --- | --- | --- | --- | --- | --- | --- |
|  | 1 | 2 | 3 | 4 | 5 | 6 | 7 | 8 | 9 | 10 |
| COSSO(GCV) | 100 | 100 | 100 | 100 | 14 | 11 | 18 | 15 | 11 | 13 |
| COSSO(5CV) | 100 | 94 | 100 | 100 | 1 | 1 | 3 | 2 | 4 | 2 |
| MARS | 100 | 100 | 100 | 100 | 35 | 35 | 34 | 39 | 28 | 35 |

$$+ 1.5g_1(x^{(5)}) + 1.5g_2(x^{(6)}) + 1.5g_3(x^{(7)}) + 1.5g_4(x^{(8)})$$
$$+ 2g_1(x^{(9)}) + 2g_2(x^{(10)}) + 2g_3(x^{(11)}) + 2g_4(x^{(12)}).$$

Therefore there are 48 uninformative variables. Let $n = 500$. The variance of the normal noise is 0.5184, to give a signal-to-noise ratio of 3:1 in the uniform case. For comparison, in the uniform setting $\text{var}\{g_1(X^{(1)})\} = 0.08$, $\text{var}\{g_2(X^{(2)})\} = 0.09$, $\text{var}\{g_3(X^{(3)})\} = 0.21$ and $\text{var}\{g_4(X^{(4)})\} = 0.26$. Both COSSO and MARS are run 100 times with additive models. The results are summarized in Table 5. We see that the two COSSO procedures outperform MARS, with the COSSO with fivefold cross-validation doing slightly better than the COSSO with GCV.

TABLE 3
*Frequency of the size of the models in the uniform setting*

|  | Model size | | | | | | | | |
| --- | --- | --- | --- | --- | --- | --- | --- | --- | --- |
|  | 3 | 4 | 5 | 6 | 7 | 8 | 9 | 10 | Mean |
| COSSO(GCV) | 0 | 57 | 17 | 18 | 5 | 2 | 0 | 1 | 4.82 |
| COSSO(5CV) | 6 | 84 | 7 | 3 | 0 | 0 | 0 | 0 | 4.07 |
| MARS | 0 | 4 | 24 | 40 | 26 | 6 | 0 | 0 | 6.06 |

TABLE 4
*Mean and standard deviation of the model sizes*

|  | Comp. symm. | | AR(1) | | |
| --- | --- | --- | --- | --- | --- |
|  | $t=1$ | $t=3$ | $\rho=-0.5$ | $\rho=0$ | $\rho=0.5$ |
| COSSO(GCV) | 4.8 (1.2) | 4.8 (1.5) | 4.7 (1.2) | 4.8 (1.3) | 4.6 (1.2) |
| COSSO(5CV) | 4.1 (1.2) | 4.4 (1.9) | 4.1 (1.2) | 4.0 (1.0) | 3.8 (0.9) |
| MARS | 6.3 (0.9) | 6.2 (0.9) | 6.1 (1.0) | 6.1 (0.8) | 5.9 (0.8) |

EXAMPLE 3. We consider a 10-dimensional regression problem with several two-way interactions:

$$f(x) = g_1(x^{(1)}) + g_2(x^{(2)}) + g_3(x^{(3)}) + g_4(x^{(4)})$$
$$+ g_1(x^{(3)}x^{(4)}) + g_2\left(\frac{x^{(1)} + x^{(3)}}{2}\right) + g_3(x^{(1)}x^{(2)}).$$

We consider the uniform setting and set the noise to be normal with standard deviation 0.2546, to give a signal-to-noise ratio of 3:1. The average integrated squared errors are given in Table 6 for sample sizes $n = 100, 200, 400$. Both the COSSO and MARS are run with the two-way interaction model. We follow the advice in [8] to set the cost for each basis function optimization to be 3 in the MARS for two-way interaction models.

There are 55 function components in the COSSO. The COSSO does not do well when $n = 100$. It seems that there are too many function components for the COSSO to select from with 100 data points. MARS does not suffer from a small sample size so much as the COSSO. Part of the reason is that the MARS algorithm introduces a certain hierarchical order of the terms being searched from: only after a univariate basis function is included in the model will the product of other terms with it become a candidate for inclusion in later steps. In contrast, the COSSO selects from all the function components and does not distinguish between main effects and interaction terms. Therefore the COSSO does not assume any hierarchical structure, and may not be efficient when the true model is hierarchical and the sample size is small. However, as the sample size increases, the COSSO procedures catch up quickly. Their performance is comparable to MARS when $n = 200$ and better than MARS when $n = 400$.

TABLE 5
*Average ISE (unit $10^{-3}$) and model sizes with their standard errors*

|  |  | Comp. symm. | | AR(1) | |
| --- | --- | --- | --- | --- | --- |
|  |  | $t = 0$ | $t = 1$ | $\rho = 0.5$ | $\rho = -0.5$ |
| ISE | COSSO(GCV) | 201 (4) | 178 (5) | 199 (6) | 183 (5) |
|  | COSSO(5CV) | 144 (4) | 162 (5) | 153 (4) | 149 (5) |
|  | MARS | 353 (7) | 302 (7) | 286 (6) | 280 (5) |
| Model sizes | COSSO(GCV) | 18.0 (4.1) | 18.0 (4.1) | 19.0 (5.1) | 18.0 (4.3) |
|  | COSSO(5CV) | 12.0 (0.2) | 11.7 (1.4) | 12.1 (1.4) | 11.9 (1.0) |
|  | MARS | 35.2 (2.3) | 36.1 (2.1) | 35.2 (2.5) | 35.9 (2.4) |

TABLE 6
*Average integrated squared errors for Example 3*

|  | $n = 100$ | $n = 200$ | $n = 400$ |
| --- | --- | --- | --- |
| COSSO(GCV) | 0.358 (0.009) | 0.100 (0.003) | 0.045 (0.001) |
| COSSO(5CV) | 0.378 (0.005) | 0.094 (0.004) | 0.043 (0.001) |
| MARS | 0.239 (0.008) | 0.109 (0.003) | 0.084 (0.001) |

In the above examples we see that in general the COSSO with fivefold cross-validation tends to do better than the COSSO with GCV. We therefore recommend the use of fivefold cross-validation with the COSSO unless the computation time is a crucial factor. In the following examples the COSSO is tuned with fivefold cross-validation.

EXAMPLE 4 (The circuit example). This is an example from [8]. Of interest is the dependence of the impedance $Z$ of a circuit and phase shift $\phi$ on components in the circuit. The true dependence is described by

$$Z = [R^2 + \{\omega L - 1/(\omega C)\}^2]^{1/2},$$
$$\phi = \tan^{-1}\left\{\frac{\omega L - 1/(\omega C)}{R}\right\}.$$

The input variables are uniform in the range $0 \leq R \leq 100$, $40\pi \leq \omega \leq 560\pi$, $0 \leq L \leq 1$, $1 \leq C \leq 11$, and the noise is normal with the standard deviation set to give a signal-to-noise level 3:1.

This is a relatively small problem with $d = 4$. All orders of interactions are present. Friedman [8] applied MARS with the additive model, the two-way interaction model and the saturated model to this example, and found that the performance of the two-way interaction model was the best. We scale the input region to $[0,1]^4$ and apply the COSSO with fivefold cross-validation. With the small dimension, it is possible to apply the COSSO with the saturated model, which has $2^4 - 1 = 15$ function components. However, it turns out that the two-way interaction COSSO does slightly better than the saturated model. We compare the integrated squared error of the two-way interaction COSSO and that of the two-way interaction MARS in Table 7. It turns out that the COSSO performs much better than MARS. We note that in Table 7 the numbers reported for MARS when $n = 200$ are calculated after excluding one large extreme outlier in ISE, and those for MARS when $n = 400$ are calculated after excluding three large extreme outliers in ISE.

**8. Real examples.** We now apply the COSSO to three real datasets: Ozone data, Boston housing data and Tecator data. The first two datasets are available from the R library "mlbench." In the Ozone data, the daily maximum one-hour-average ozone reading and eight meteorological variables were recorded in the Los Angeles basin for 330 days in 1976. The Boston housing data concerns housing values in the suburbs of Boston. There are 12 input variables. The sample size is 506. The Tecator data is available from the datasets archive of StatLib at [lib.stat.cmu.edu/datasets/](lib.stat.cmu.edu/datasets/). The data was recorded on a Tecator Infratec Food and Feed Analyzer working in the wavelength range 850–1050 nm by the Near Infrared Transmission (NIT) principle. Each sample contains finely chopped pure meat with different moisture, fat and protein content. The input vector consists of a 100-channel spectrum of absorbances. The absorbance is $-\log_{10}$ of the transmittance measured by the spectrometer. As suggested in the document, we use the first 13 principal components to predict the fat content. The total sample size is 215.

TABLE 7
*Average integrated squared error for estimating the impedance $Z$ (in units of $10^3$) and the phase shift $\phi$ (in units of $10^{-3}$)*

|  |  | $n=100$ | $n=200$ | $n=400$ |
|---|---|---|---|---|
| for $Z$ | COSSO | 1.91 (0.12) | 0.85 (0.05) | 0.51 (0.03) |
|  | MARS | 5.57 (0.41) | 2.47 (0.16) | 1.37 (0.08) |
| for $\phi$ | COSSO | 12.98 (0.36) | 7.96 (0.20) | 5.36 (0.10) |
|  | MARS | 20.59 (0.96) | 12.60 (0.71) | 8.19 (0.14) |

TABLE 8
*Estimated prediction squared errors and their standard errors*

|  | Ozone | Boston | Tecator |
|---|---|---|---|
| COSSO | 16.04 (0.06) | 9.89 (0.08) | 0.92 (0.02) |
| MARS | 18.24 (0.45) | 14.31 (0.34) | 4.99 (1.07) |

We apply the COSSO and MARS on these datasets, and estimate the prediction squared errors $E[\{Y - \hat{f}(X)\}^2]$ by tenfold cross-validation. We select the tuning parameter by fivefold CV within the training set. The estimate obtained is then evaluated on the test set. We do this tenfold cross-validation five times and then average. For all three datasets, both the COSSO and MARS find the two-way interaction model has better prediction accuracy than the additive model. Therefore we choose to apply the two-way interaction model and the results are given in Table 8. We can see the COSSO does considerably better than MARS.

The COSSO algorithm can be naturally extended to tackle the nonparametric classification problem where the response $y$ is either 0 or 1. In this case we replace the square loss in (6) by the logistic regression loss function. The first term in (6) is then replaced by $\ell(f) = 1/n \sum_{i=1}^{n}[-y_i f(x_i) + \log(1 +$ $e^{f(x_i)})]$, where $f$ is the logit function. We can then solve the optimization problem by applying the quadratic approximation to $\ell(f)$ iteratively. This leads to the iteratively reweighted least squares (IRLS) procedure, which is equivalent to a Newton–Raphson algorithm. Using this approach, we can solve our optimization problem by iterative application of the COSSO algorithm, within an IRLS loop. We illustrate the performance of this procedure by comparing it with a number of machine learning algorithms on several high-dimensional real datasets.

van Gestel et al. [18] conducted a benchmark study comparing a number of commonly used machine learning techniques including the support vector machine (SVM), least squares SVM (LS-SVM), linear discriminant analysis (LDA), quadratic discriminant analysis (QDA), logistic regression (Logit), the decision tree algorithm C4.5, Holte's one-rule classifier (oneR), instance-based learners (IB) and the Naive Bayes method. There are five binary classification datasets with continuous predictors, and we test the performance of the COSSO on these datasets. The datasets are the BUPA Liver Disorder data, the Johns Hopkins University Ionosphere data, the PIMA Indian Diabetes data, the Sonar, Mines vs. Rocks data and the Wisconsin Breast Cancer data. The basic features of the datasets and the performance of different algorithms are summarized in Table 9. Due to the high dimension of these problems, we only consider the COSSO additive model.

Following [18], for each dataset we randomly select 2/3 of the data for training and tuning, and test on the remaining 1/3 of the data. We do this randomization ten times and report the average test set performance and sample standard deviation for the COSSO. The results for the other algorithms are taken from [18]. van Gestel et al. [18] included six types of LS-SVM's and found that the LS-SVM with the radial basis function (RBF) kernel performs best overall. To save space we only include the LS-SVM with RBF kernel and linear kernel in Table 9. van Gestel et al. [18] included two instance-based learners (IB1 and IB10 in [18]). We combine them and report the better performance of the two. The same is done for the two Naive Bayes methods considered in [18]. The best average test set performance is denoted in boldface for each dataset in Table 9. We can see the COSSO gives very competitive performance on these benchmark datasets.

**9. Discussion.** The difference between the COSSO and the common smoothing spline in smoothing spline ANOVA mirrors that between the LASSO and ridge regression in linear models. Compared with other model selection algorithms based on greedy search, the COSSO optimizes a global criterion and provides a shrinkage estimate. We have shown that the COSSO has attractive properties for model selection and estimation.

One future research topic is statistical inference based on the COSSO. Traditionally inference after model selection is based on the selected models, resulting in biased inference. Shen, Huang and Ye [12] proposed a method to make approximately unbiased inference. It is of interest to see if their method can be adapted to our problem to give unbiased inference.

TABLE 9
*Comparison of the test set performance of the COSSO with the performance of SVM, LS-SVM, LDA, QDA, Logit, C4.5, oneR, IB, Naive Bayes and Majority Rule*

|                  | BUPA          | Ionosphere    | Pima Indian   | Sonar MR      | Wisc. BC      |
|------------------|---------------|---------------|---------------|---------------|---------------|
| $n$              | 345           | 351           | 768           | 208           | 683           |
| $d$              | 6             | 33            | 8             | 60            | 9             |
| COSSO            | **71.1** (3.5)| 91.1 (3.7)    | **77.3** (2.2)| **79.0** (4.5)| 97.0 (0.8)    |
| SVM (linear)     | 67.7 (2.6)    | 87.1 (3.4)    | 77.0 (2.4)    | 74.1 (4.2)    | 96.3 (1.0)    |
| SVM (RBF)        | 70.4 (3.2)    | 95.4 (1.7)    | **77.3** (2.2)| 75.0 (6.6)    | 96.4 (1.0)    |
| LS-SVM (linear)  | 65.6 (3.2)    | 87.9 (2.0)    | 76.8 (1.8)    | 72.6 (3.7)    | 95.8 (1.0)    |
| LS-SVM (RBF)     | 70.2 (4.1)    | **96.0** (2.1)| 76.8 (1.7)    | 73.1 (4.2)    | 96.4 (1.0)    |
| LDA              | 65.4 (3.2)    | 87.1 (2.3)    | 76.7 (2.0)    | 67.9 (4.9)    | 95.6 (1.1)    |
| QDA              | 62.2 (3.6)    | 90.6 (2.2)    | 74.2 (3.3)    | 53.6 (7.4)    | 94.5 (0.6)    |
| Logit            | 66.3 (3.1)    | 86.2 (3.5)    | 77.2 (1.8)    | 68.4 (5.2)    | 96.1 (1.0)    |
| C4.5             | 63.1 (3.8)    | 90.6 (2.2)    | 73.5 (3.0)    | 72.1 (2.5)    | 94.7 (1.0)    |
| oneR             | 56.3 (4.4)    | 83.6 (4.8)    | 71.3 (2.7)    | 62.6 (5.5)    | 91.8 (1.4)    |
| IB               | 61.3 (6.2)    | 87.2 (2.8)    | 73.6 (2.4)    | 77.7 (4.4)    | 96.4 (1.2)    |
| Naive Bayes      | 63.7 (4.5)    | 92.1 (2.5)    | 75.5 (1.7)    | 71.6 (3.5)    | **97.1** (0.9)|
| Majority Rule    | 56.5 (3.1)    | 64.4 (2.9)    | 66.8 (2.1)    | 54.4 (4.7)    | 66.2 (2.4)    |

## APPENDIX

**Proofs.**

PROOF OF THEOREM 1. Denote the functional to be minimized in (6) by $A(f)$; then $A(f)$ is convex and continuous. Without loss of generality, we assume $\tau = 1$.

By (7) we have that $J(f) \geq \|f\|$, for any $f \in \mathcal{F}_1$. Let $R_{\mathcal{F}_1}$ be the reproducing kernel of $\mathcal{F}_1$ and let $\langle \cdot, \cdot \rangle_{\mathcal{F}_1}$ be the inner product in $\mathcal{F}_1$. Denote $a = \max_{i=1}^{n} R_{\mathcal{F}_1}^{1/2}(x_i, x_i)$. By the definition of a reproducing kernel, we have for any $f \in \mathcal{F}_1$ and $i = 1, \ldots, n$,

$$
\begin{aligned}
(A.1) \quad |f(x_i)| &= |\langle f(\cdot), R_{\mathcal{F}_1}(x_i, \cdot) \rangle_{\mathcal{F}_1}| \leq \|f\| \langle R_{\mathcal{F}_1}(x_i, \cdot), R_{\mathcal{F}_1}(x_i, \cdot) \rangle_{\mathcal{F}_1}^{1/2} \\
&= \|f\| R_{\mathcal{F}_1}^{1/2}(x_i, x_i) \leq a\|f\| \leq aJ(f).
\end{aligned}
$$

Denote $\rho = \max_{i=1}^{n}(y_i^2 + |y_i| + 1)$. Consider the set

$D = \{f \in \mathcal{F} : f = b + f_1,$

with $b \in \{1\}, f_1 \in \mathcal{F}_1, J(f) \leq \rho, |b| \leq \rho^{1/2} + (a+1)\rho\}.$

Then $D$ is a closed, convex and bounded set. Therefore by Theorem 4 of [14], page 162, there exists a minimizer of (6) in $D$. Denote the minimizer by $\bar{f}$. Then $A(\bar{f}) \leq A(0) < \rho$.

On the other hand, for any $f \in \mathcal{F}$ with $J(f) > \rho$, clearly $A(f) \geq J(f) > \rho$; for any $f \in \mathcal{F}$ with $J(f) \leq \rho$, $f = b + f_1$, $b \in \{1\}$, $f_1 \in \mathcal{F}$ and $|b| > \rho^{1/2} + (a+1)\rho$, we use (A.1) to get that, for any $i = 1, \ldots, n$,

$$|b + f_1(x_i) - y_i| > [\rho^{1/2} + (a+1)\rho] - a\rho - \rho = \rho^{1/2}.$$

Therefore $A(f) > \rho$. For any $f \notin D$, $A(f) > A(\bar{f})$, that is, $\bar{f}$ is a minimizer of (6) in $\mathcal{F}$. □

PROOF OF THEOREM 2. For any $f \in \mathcal{F}$, we can write $f(x) = c + f_1(x^{(1)}) + \cdots + f_d(x^{(d)}) = c + g(x)$, such that $\sum_{i=1}^n f_j(x_i^{(j)}) = 0$, $j = 1, \ldots, d$. Similarly, write $f_0(x) = c_0 + f_{01}(x^{(1)}) + \cdots + f_{0d}(x^{(d)}) = c_0 + g_0(x)$ and $\hat{f}(x) = \hat{c} + \hat{f}_1(x^{(1)}) + \cdots + \hat{f}_d(x^{(d)}) = \hat{c} + \hat{g}(x)$. By construction $\sum_{i=1}^n \{g_0(x_i) - g(x_i)\} = 0$, and we can write (6) as

$$(c_0 - c)^2 + \frac{2}{n}(c_0 - c)\sum_{i=1}^n \varepsilon_i + \frac{1}{n}\sum_{i=1}^n \{g_0(x_i) + \varepsilon_i - g(x_i)\}^2 + \tau_n^2 J(g).$$

Therefore, the minimizing $\hat{c}$ must minimize $(c_0 - c)^2 + 2/n(c_0 - c)\sum_{i=1}^n \varepsilon_i$. That is, $\hat{c} = c_0 + 1/n\sum_{i=1}^n \varepsilon_i$. Therefore $(\hat{c} - c_0)^2$ converges with rate $n^{-1}$. On the other hand, $\hat{g}$ must minimize

$$\frac{1}{n}\sum_{i=1}^n \{g_0(x_i) + \varepsilon_i - g(x_i)\}^2 + \tau_n^2 J(g).$$

Let $\mathcal{G} = \{g \in \mathcal{F} : g(x) = f_1(x^{(1)}) + \cdots + f_d(x^{(d)}), \text{ with } \sum_{i=1}^n f_j(x_i^{(j)}) = 0, j = 1, \ldots, d\}$. Then $g_0 \in \mathcal{G}$, $\hat{g} \in \mathcal{G}$. The conclusion of Theorem 2 then follows from Theorem 10.2 of [17] and the following lemma. □

LEMMA A.1. *Let $H_\infty(\delta, \mathcal{G})$ be the $\delta$-entropy of $\mathcal{G}$ for the supremum norm. Then*

$$H_\infty(\delta, \{g \in \mathcal{G} : J(g) \leq 1\}) \leq A d^{(\ell+1)/\ell} \delta^{-1/\ell},$$

*for all $\delta > 0$, $n \geq 1$ and some $A > 0$ not depending on $\delta$, $n$ or $d$.*

PROOF. Define $\mathcal{G}^j$ as the set of univariate functions of $x^{(j)}$

$$\mathcal{G}^j = \left\{f_j \in S_\ell : J(f_j) \leq 1, \sum_{i=1}^n f_j(x_i^{(j)}) = 0\right\},$$

COMPONENT SELECTION AND SMOOTHING 23where $S_\ell$ is the $\ell$th-order Sobolev space. Then from (3), any $h \in \mathcal{G}^j$ satisfies

$$(A.2) \qquad \sum_{\nu=0}^{\ell-2}[h^{(\nu)}(1) - h^{(\nu)}(0)]^2 + \int_0^1 \{h^{(\ell)}(t)\}^2 \, dt \leq 1.$$

We first show that for any $h \in \mathcal{G}^j$, we have $|h|_\infty \equiv \{\sup_{s \in [0,1]} |h(s)|\} \leq 1$. For any fixed $h \in \mathcal{G}^j$, define $K = \{k : 0 \leq k \leq \ell - 1,$ and for any integer $q \in [0, k]$, there exists $a_q \in [0, 1]$ satisfying $h^{(q)}(a_q) = 0\}$. Since $\sum_{i=1}^n h(x_i^{(j)}) = 0$, we have $0 \in K$. Let $k_0$ be the largest number in $K$. Now we consider the two possibilities $k_0 \neq \ell - 1$ and $k_0 = \ell - 1$ separately.

If $k_0 \neq \ell - 1$, then $k_0 \leq \ell - 2$, $k_0 \in K$ and $k_0 + 1 \notin K$. By the definition of $K$, we see that $h^{(k_0)}$ is monotone and crosses the $x$-axis. So for any $s \in [0, 1]$, we have $|h^{(k_0)}(s)| \leq |h^{(k_0)}(1) - h^{(k_0)}(0)| \leq 1$. The last inequality follows from (A.2). From this and the fact that $h^{(k_0-1)}$ crosses the $x$-axis, we get $|h^{(k_0-1)}(s)| \leq 1$, $\forall s \in [0, 1]$. Continuing with this argument, we get $|h|_\infty \leq 1$.

On the other hand, if $k_0 = \ell - 1$, then $\max_s h^{(\ell-1)}(s) \geq 0$, and $\min_s h^{(\ell-1)}(s) \leq 0$. By (A.2) we get

$$1 \geq \int_0^1 \{h^{(\ell)}(t)\}^2 \, dt \geq \left\{\int_0^1 |h^{(\ell)}(t)| \, dt\right\}^2 \geq \left\{\max_s h^{(\ell-1)}(s) - \min_s h^{(\ell-1)}(s)\right\}^2.$$

Therefore $-1 \leq \min_s h^{(\ell-1)}(s) \leq 0 \leq \max_s h^{(\ell-1)}(s) \leq 1$. That is, $|h^{(\ell-1)}(s)| \leq 1$, $\forall s \in [0, 1]$. Now by the definition of $K$ and that $\ell - 1 \in K$, we know that $h^{(k)}$ crosses the $x$-axis for any integer $k \in [0, \ell-1]$. Therefore we get $|h|_\infty \leq 1$.

Therefore we have shown that $|h|_\infty \leq 1$ for any $h \in \mathcal{G}^j$. It then follows from Theorem 2.4 of [17], page 19, that

$$(A.3) \qquad H_\infty(\delta, \mathcal{G}^j) \leq A\delta^{-1/\ell}$$

for all $\delta > 0$ and $n \geq 1$, and some positive $A$ not depending on $\delta$ and $n$.

By the definition of $\mathcal{G}$ and the $\mathcal{G}^j$, we see that in terms of the supreme norm, if each $\mathcal{G}^j$, $j = 1, \ldots, d$, can be covered by $N$ balls of radius $\delta$, then the set $\{g \in \mathcal{G} : J(g) \leq 1\}$ can be covered by $N^d$ balls with radius $d\delta$. By (A.3) we get

$$H_\infty(d\delta, \{g \in \mathcal{G} : J(g) \leq 1\}) \leq Ad\delta^{-1/\ell},$$

and the conclusion of the lemma follows. $\square$

PROOF OF LEMMA 1. For any $f \in \mathcal{F}$, we can write $f = b + \sum_{\alpha=1}^p f_\alpha$ with $f_\alpha \in \mathcal{F}^\alpha$. Let the projection of $f_\alpha$ onto $\mathrm{span}\{R_\alpha(x_i, \cdot), i = 1, \ldots, n\} \subset \mathcal{F}^\alpha$ be $g_\alpha$ and its orthogonal complement be $h_\alpha$. Then $f_\alpha = g_\alpha + h_\alpha$, and $\|f_\alpha\|^2 =$ $\|g_\alpha\|^2 + \|h_\alpha\|^2$, $\alpha = 1, \ldots, p$. Since $R = 1 + \sum_{\alpha=1}^{p} R_\alpha$ is the reproducing kernel of $\mathcal{F}$, we have, making use of the orthogonal structures,

$$f(x_i) = \left\langle 1 + \sum_{\alpha=1}^{p} R_\alpha(x_i, \cdot), b + \sum_{\alpha=1}^{p} (g_\alpha + h_\alpha) \right\rangle = b + \sum_{\alpha=1}^{p} \langle R_\alpha(x_i, \cdot), g_\alpha \rangle,$$

where $\langle \cdot, \cdot \rangle$ is the inner product in $\mathcal{F}$. Therefore (6) can be written as

$$\frac{1}{n} \sum_{i=1}^{n} \left\{ y_i - b - \sum_{\alpha=1}^{p} \langle R_\alpha(x_i, \cdot), g_\alpha \rangle \right\}^2 + \tau^2 \sum_{\alpha=1}^{p} (\|g_\alpha\|^2 + \|h_\alpha\|^2)^{1/2}.$$

Therefore any minimizing $f$ satisfies $h_\alpha = 0$, $\alpha = 1, \ldots, p$, and the conclusion of the lemma follows. □

PROOF OF LEMMA 2. Denote the functional in (6) by $A(f)$ and the functional in (8) by $B(\theta, f)$. We have $\lambda_0 \theta_\alpha^{-1} \|P^\alpha f\|^2 + \lambda \theta_\alpha \geq 2\lambda_0^{1/2} \lambda^{1/2} \|P^\alpha f\| = \tau^2 \|P^\alpha f\|$, for any $\theta_\alpha \geq 0$ and $f \in \mathcal{F}$, and equality holds if and only if $\theta_\alpha = \lambda_0^{1/2} \lambda^{-1/2} \|P^\alpha f\|$. Therefore $B(\theta, f) \geq A(f)$ for any $\theta_\alpha \geq 0$ and $f \in \mathcal{F}$, and equality holds if and only if $\theta_\alpha = \lambda_0^{1/2} \lambda^{-1/2} \|P^\alpha f\|$, $\alpha = 1, \ldots, p$. The conclusion of the lemma follows. □

**Further derivations in the tensor product design case.** Now we consider the function space $\mathcal{F} = T^1 \otimes T^2$. Define $\Sigma = \{\bar{K}(x_{i,1}, x_{j,1})\}_{m \times m}$, the marginal kernel matrix corresponding to the reproducing kernel of $\bar{T}$. With a little abuse of notation, let $R_j$, $j = 1, 2$, also stand for the $n \times n$ matrix of the reproducing kernel $R_j$ evaluated at the $n$ data points, and the same for $R_{12}$. Suppose the data points are permuted appropriately; we have $R_{12} = \Sigma \otimes \Sigma$, where $\otimes$ stands for the Kronecker product of matrices. Let $\mathbf{1}_m$ be the column vector consisting of $m$ 1's. For the main effect spaces, we have $R_1 = \Sigma \otimes (\mathbf{1}_m \mathbf{1}_m^T)$ and $R_2 = (\mathbf{1}_m \mathbf{1}_m^T) \otimes \Sigma$.

Straightforward calculation gives $\Sigma \mathbf{1}_m = mt \mathbf{1}_m$, where $t = 1/(720 m^4)$. Let $\{\xi_1 = \mathbf{1}_m, \xi_2, \ldots, \xi_m\}$ be an orthonormal (with respect to the inner product $\langle \cdot, \cdot \rangle_m$ in $R^m$) eigensystem of $\Sigma$, with corresponding eigenvalues $m\eta_1, m\eta_2, \ldots, m\eta_m$, where $\eta_1 = t$, and $\eta_2 \geq \eta_3 \geq \cdots \geq \eta_m$. Then it is well known that $\eta_i \sim i^{-4}$ for $i \geq 2$. See [16]. Notice $\xi_1, \xi_2, \ldots, \xi_m$ are also the eigenvectors of $\mathbf{1}_m \mathbf{1}_m^T$, with corresponding eigenvalues being $m, 0, \ldots, 0$. Write $\xi_{\mu\nu} = \xi_\mu \otimes \xi_\nu$. It is then easy to check that $\{\xi_{\mu\nu} : \mu, \nu = 1, \ldots, m\}$ form an eigensystem of $R_1$, $R_2$ and $R_{12}$. The eigenvalues of $R_1$, $R_2$ and $R_{12}$ are, respectively,

$$\begin{aligned}
r_{1,\mu 1} &= n\eta_\mu; & r_{1,\mu\nu} &= 0 & &\text{for } \mu \geq 1, \nu \geq 2; \\
r_{2,1\nu} &= n\eta_\nu; & r_{2,\mu\nu} &= 0 & &\text{for } \mu \geq 2, \nu \geq 1; \\
r_{12,\mu\nu} &= n\eta_\mu \eta_\nu & & & &\text{for } \mu \geq 1, \nu \geq 1.
\end{aligned}$$

It is clear that $\{\xi_{\mu\nu}: \mu,\nu = 1,\ldots,m\}$ is also an orthonormal basis in $R^n$ with respect to the inner product $\langle \cdot,\cdot\rangle_n$. Consider the vector of length $n$ of function values at the sample points: $f = (f(x_{k,1}, x_{\ell,2}): k,\ell = 1,\ldots,m)^{\mathrm{T}}$. Let $O$ be the $n \times n$ matrix with columns being the vectors $\xi_{\mu\nu}$, $\mu,\nu = 1,\ldots,m$. Then $O^{\mathrm{T}}O = nI$. Denote $a = (a_{\mu\nu}: \mu,\nu = 1,\ldots,m)^{\mathrm{T}} = (1/n)O^{\mathrm{T}}f$ and $z = (z_{\mu\nu}: \mu,\nu = 1,\ldots,m)^{\mathrm{T}} = (1/n)O^{\mathrm{T}}y$. That is,

$$a_{\mu\nu} = \langle f, \xi_{\mu\nu}\rangle_n, \qquad z_{\mu\nu} = \langle y, \xi_{\mu\nu}\rangle_n;$$

then $f \in R^n$ can be expanded in terms of the orthonormal basis,

$$f = \sum_{\mu,\nu} a_{\mu\nu}\xi_{\mu\nu} = f_0 + f_1 + f_2 + f_{12},$$

where $f_0 = a_{11}\xi_{11}$, $f_1 = \sum_{\mu=2}^{m} a_{\mu 1}\xi_{\mu 1}$, $f_2 = \sum_{\nu=2}^{m} a_{1\nu}\xi_{1\nu}$ and $f_{12} = \sum_{\mu=2}^{m}\sum_{\nu=2}^{m} a_{\mu\nu}\xi_{\mu\nu}$. Then all the components $f_0$, $f_1$, $f_2$ and $f_{12}$ are orthogonal in $R^n$. Furthermore, we have $z_{\mu\nu} = a_{\mu\nu} + \delta_{\mu\nu}$, where $\delta_{\mu\nu} \sim N(0,\sigma^2/n)$, for $\mu \geq 1, \nu \geq 1$.

Now let us consider the COSSO estimate (11),

$$\frac{1}{n}(y - R_\theta c - b\mathbf{1}_n)^{\mathrm{T}}(y - R_\theta c - b\mathbf{1}_n) + c^{\mathrm{T}}R_\theta c + \lambda\sum_{\alpha=1}^{p}\theta_\alpha \quad \text{subject to} \quad \theta_\alpha \geq 0,$$

$$\alpha = 1,\ldots,p,$$

where $R_\theta = \sum_{\alpha=1}^{p}\theta_\alpha R_\alpha$. Let $s = O^{\mathrm{T}}c$, $D_\alpha = (1/n^2)O^{\mathrm{T}}R_\alpha O$. Then $D_\alpha$ is a diagonal matrix with diagonal elements $r_{\alpha,\mu\nu}/n$, $\alpha = 1, 2$ or $12$. The COSSO problem can be written as

$$(z - D_\theta s - (b,0,\ldots,0)^{\mathrm{T}})^{\mathrm{T}}(z - D_\theta s - (b,0,\ldots,0)^{\mathrm{T}}) + s^{\mathrm{T}}D_\theta s + \lambda\sum_{\alpha=1}^{p}\theta_\alpha,$$

where $D_\theta = \sum_{\alpha=1}^{p}\theta_\alpha D_\alpha$. It can then be shown by straightforward calculation that, for the minimizing $s$, $b$ and $\theta$, $\hat{s}_{11} = 0$, $\hat{b} = z_{11}$ and $s$ and $\theta$ minimize

$$\sum_{\mu \geq 2}[\{z_{\mu 1} - \eta_\mu(\theta_1 + \theta_{12}t)s_{\mu 1}\}^2 + \eta_\mu(\theta_1 + \theta_{12}t)s_{\mu 1}^2]$$

$$+ \sum_{\nu \geq 2}[\{z_{1\nu} - \eta_\nu(\theta_2 + \theta_{12}t)s_{1\nu}\}^2 + \eta_\nu(\theta_2 + \theta_{12}t)s_{1\nu}^2]$$

$$+ \sum_{\mu \geq 2, \nu \geq 2}[(z_{\mu\nu} - \theta_{12}\eta_\mu\eta_\nu s_{\mu\nu})^2 + \theta_{12}\eta_\mu\eta_\nu s_{\mu\nu}^2] + \lambda(\theta_1 + \theta_2 + \theta_{12}).$$

Therefore, at the minimum we have

$$\hat{s}_{\mu 1} = \{1 + \eta_\mu(\theta_1 + \theta_{12}t)\}^{-1}z_{\mu 1}, \qquad \mu \geq 2;$$
$$\hat{s}_{1\nu} = \{1 + \eta_\nu(\theta_2 + \theta_{12}t)\}^{-1}z_{1\nu}, \qquad \nu \geq 2;$$
$$\hat{s}_{\mu\nu} = (1 + \eta_\mu\eta_\nu\theta_{12})^{-1}z_{\mu\nu}, \qquad \mu \geq 2, \nu \geq 2;$$

and the $\theta$'s minimize

$$A(\theta_1, \theta_2, \theta_{12}) = \sum_{\mu \geq 2} z_{\mu 1}^2 (1 + \eta_\mu \theta_1 + \eta_\mu \theta_{12} t)^{-1} + \sum_{\nu \geq 2} z_{1\nu}^2 (1 + \eta_\nu \theta_2 + \eta_\nu \theta_{12} t)^{-1}$$
$$+ \sum_{\mu \geq 2, \nu \geq 2} z_{\mu\nu}^2 (1 + \theta_{12} \eta_\mu \eta_\nu)^{-1} + \lambda(\theta_1 + \theta_2 + \theta_{12}),$$

subject to $\theta_1 \geq 0, \theta_2 \geq 0, \theta_{12} \geq 0$. A calculation similar to that in Section 4 then shows that, when $\lambda$ is appropriately chosen, the COSSO selects the correct components with probability tending to 1.

**Acknowledgment.** The authors wish to thank Grace Wahba for helpful comments.

## REFERENCES

[1] BREIMAN, L. (1995). Better subset selection using the nonnegative garrote. *Technometrics* **37** 373–384. MR1365720
[2] CHEN, S., DONOHO, D. and SAUNDERS, M. (1998). Atomic decomposition by basis pursuit. *SIAM J. Sci. Comput.* **20** 33–61. MR1639094
[3] CHEN, Z. (1993). Fitting multivariate regression functions by interaction spline models. *J. Roy. Statist. Soc. Ser. B* **55** 473–491. MR1224411
[4] CRAVEN, P. and WAHBA, G. (1979). Smoothing noisy data with spline functions. *Numer. Math.* **31** 377–403. MR0516581
[5] EFRON, B., HASTIE, T., JOHNSTONE, I. and TIBSHIRANI, R. (2004). Least angle regression (with discussion). *Ann. Statist.* **32** 407–499. MR2060166
[6] FAN, J. and LI, R. (2001). Variable selection via nonconcave penalized likelihood and its oracle properties. *J. Amer. Statist. Assoc.* **96** 1348–1360. MR1946581
[7] FRANK, I. E. and FRIEDMAN, J. H. (1993). A statistical view of some chemometrics regression tools. *Technometrics* **35** 109–148.
[8] FRIEDMAN, J. H. (1991). Multivariate adaptive regression splines (with discussion). *Ann. Statist.* **19** 1–141. MR1091842
[9] GU, C. (1992). Diagnostics for nonparametric regression models with additive terms. *J. Amer. Statist. Assoc.* **87** 1051–1058.
[10] GU, C. (2002). *Smoothing Spline ANOVA Models*. Springer, Berlin. MR1876599
[11] GUNN, S. R. and KANDOLA, J. S. (2002). Structural modeling with sparse kernels. *Machine Learning* **48** 137–163.
[12] SHEN, X., HUANG, H. and YE, J. (2004). Inference after model selection. *J. Amer. Statist. Assoc.* **99** 751–762. MR2090908
[13] SHEN, X. and YE, J. (2002). Adaptive model selection. *J. Amer. Statist. Assoc.* **97** 210–221. MR1947281
[14] TAPIA, R. and THOMPSON, J. (1978). *Nonparametric Probability Density Estimation*. Johns Hopkins Univ. Press, Baltimore. MR0502724
[15] TIBSHIRANI, R. J. (1996). Regression shrinkage and selection via the lasso. *J. Roy. Statist. Soc. Ser. B* **58** 267–288. MR1379242
[16] UTRERAS, F. (1983). Natural spline functions: Their associated eigenvalue problem. *Numer. Math.* **42** 107–117. MR0716477
[17] VAN DE GEER, S. (2000). *Empirical Processes in M-Estimation*. Cambridge Univ. Press.

[18] van Gestel, T., Suykens, J. A. K., Baesens, B., Viaene, S., Vanthienen, J., Dedene, G., de Moor, B. and Vandewalle, J. (2004). Benchmarking least squares support vector machine classifiers. *Machine Learning* **54** 5–32.
[19] Wahba, G. (1990). *Spline Models for Observational Data*. SIAM, Philadelphia. MR1045442
[20] Wahba, G., Wang, Y., Gu, C., Klein, R. and Klein, B. (1995). Smoothing spline ANOVA for exponential families, with application to the Wisconsin Epidemiological Study of Diabetic Retinopathy. *Ann. Statist.* **23** 1865–1895. MR1389856
[21] Yau, P., Kohn, R. and Wood, S. (2003). Bayesian variable selection and model averaging in high-dimensional multinomial nonparametric regression. *J. Comput. Graph. Statist.* **12** 23–54. MR1965210
[22] Ye, J. (1998). On measuring and correcting the effects of data mining and model selection. *J. Amer. Statist. Assoc.* **93** 120–131. MR1614596
[23] Zhang, H. H., Wahba, G., Lin, Y., Voelker, M., Ferris, M., Klein, R. and Klein, B. (2004). Variable selection and model building via likelihood basis pursuit. *J. Amer. Statist. Assoc.* **99** 659–672. MR2090901

Department of Statistics
University of Wisconsin–Madison
Madison, Wisconsin 53706
USA
E-mail: yilin@stat.wisc.edu
URL: http://www.stat.wisc.edu/~yilin

Department of Statistics
North Carolina State University
Raleigh, North Carolina 27695-8203
USA
E-mail: hzhang2@stat.ncsu.edu
URL: www4.stat.ncsu.edu/~hzhang